\documentclass[12pt,a4paper]{amsart}

\usepackage{latexsym} 
 
\usepackage{graphicx}
\usepackage{epsfig}
\usepackage{mathrsfs}
\usepackage{amssymb}
\usepackage{amsthm}
\usepackage{amsfonts}
\usepackage{amsmath}
\usepackage{amstext}
\usepackage{amscd}
\usepackage{tikz}
\usepackage{epic}
\usepackage{eepic}
\usepackage{pstricks,pst-plot}

\numberwithin{equation}{section}

\theoremstyle{plain}
\newtheorem{thm}{Theorem}[section] 
\newtheorem{prop}[thm]{Proposition}
\newtheorem{cor}[thm]{Corollary}
\newtheorem{lem}[thm]{Lemma}
\newtheorem{ass}[thm]{Assertion}

\newtheorem{theorem*}{Theorem}[]

\theoremstyle{definition}
\newtheorem{defn}[thm]{Definition}

\newtheorem{notation}{Notation}
\newtheorem{example}[thm]{Example}

\newtheorem{ques}{Question}

\theoremstyle{remark}
\newtheorem{rem}[thm]{Remark}
\theoremstyle{property}
\newtheorem{pro}[thm]{Property}

\newcommand{\N}{\mathbb{N}}
\newcommand{\R}{\mathbb{R}}

\def\accentclass@{7}
\def\makeacc@#1#2{\def#1{\mathaccent"\accentclass@#2 }}
\makeacc@\cir{017}

\title[Stabilisation of Geometric Directional Bundle]
{Stabilisation of Geometric Directional Bundle \\ 
for a subanalytic set}

\author{Satoshi Koike and Laurentiu Paunescu}

\address{Department of Mathematics, Hyogo University of Teacher Education,
942-1 Shimokume, Kato, Hyogo 673-1494, Japan}
\email{koike@hyogo-u.ac.jp} 
\address{School of Mathematics and Statistics, University of Sydney, Sydney, 
NSW, 2006, Australia}
\email{laurentiu.paunescu@sydney.edu.au}


\date{\today}

\begin{document}

\thanks{This research is partially supported by JSPS KAKENHI Grant
Number JP20K03611}


\begin{abstract}
In the previous paper \cite{kp6} we have introduced the notion of geometric 
directional bundle of a singular space, in order to introduce global 
bi-Lipschitz invariants.
Then we have posed the question of whether or not the geometric directional 
bundle is stabilised  as an operation acting on singular spaces.
In this paper we give a positive answer in the case where the singular 
spaces are subanalytic sets, thus providing a new invariant associated with the subanalytic sets.
\end{abstract}

\maketitle
\section{Introduction.}\label{introduction}

In a series of papers \cite{kp1, kp2, kp3, kp4, kp5}, in order to introduce several local bi-Lipschitz invariants,
we investigated  several local directional properties.
Subsequently, in order to introduce global bi-Lipschitz invariants, 
we have globalised the local directional properties in \cite{kp6}.
In addition, to further study the global directional properties, we have  
introduced the notion of the geometric directional bundle of a singular space.
We denote by $\mathcal{GD}(A)$ the geometric directional bundle of $A$ 
over $\{ 0 \},$ for $A \subset \R^n$ with $0 \in \overline{A}$
(see Definition \ref{gdirectionset} in \S \ref{GDB} for the definition 
of geometric directional bundle).
Then we have posed the following question in \cite{kp6}.

\begin{ques}\label{question11}(Question 1 in \cite{kp6}).
Is the operator $\mathcal{GD}$ stabilised?
Namely, does there exist a natural number $m \in \N$, not depending on $A$ 
but on $n$, such that
$$
\mathcal{GD}^m(A) = \mathcal{GD}^{m+1}(A) = \mathcal{GD}^{m+2}(A) = \cdots \ \ ?
$$
\end{ques}
\noindent
Throughout this paper, let us denote by $\N$ the set of natural numbers 
in the sense of positive integers.

In the case where $n = 2$, we have obtained a positive answer to the above 
question.

\begin{prop}\label{2stabilisation}$($Proposition 5.1 in \cite{kp6}$)$.
Let $A \subset \R^2$ such that $0 \in \overline{A}$.
Then we have 
$$
\mathcal{GD}^3 (A) = \mathcal{GD}^4 (A) = \mathcal{GD}^5 (A) = \cdots \ .
$$ 
In other words, the operation $\mathcal{GD}$ is stabilised at degree $3$.
\end{prop}

In fact, in the  Example 5.2 in \cite{kp6} we have constructed a singular set 
$A \subset \R^2$ with $0 \in \overline{A}$ such that 
the operation $\mathcal{GD}$ is stabilised at degree $3$ but not 
degree $2$ for  $A$. 
Therefore we can see that $3$ is the best degree of stabilisation 
for the operation $\mathcal{GD},$ for a general $A \subset \R^2$ with 
$0 \in \overline{A}$.

In this paper we give a positive answer to  the Question \ref{question11} 
in the case where $A \subset \R^n$ is a subanalytic set.
(See H. Hironaka \cite{hironaka} for the definition of a subanalytic set.)
Our main result is the following.

\begin{thm}\label{weakconj}
Let $A \subset \R^n$ be a subanalytic set such that 
$0 \in \overline{A}$. 
Then we have 
$$
\mathcal{GD}^{n-1}(A) = \mathcal{GD}^{n}(A) = \mathcal{GD}^{n+1}(A) = \cdots  
$$
if $n \ge 2$, and $\mathcal{GD}$ is stabilised at degree $1$ if $n = 1$.
\end{thm}

In \S3 we shall prove Theorem \ref{weakconj}, in fact it is a straightforward consequence of our Proposition \ref{main}.
Concerning this theorem, we may ask whether we can strengthen the statement or not.
Namely, we pose the following question.

\begin{ques}\label{SC}
Does there exist a natural number $m \in \N$ such that 
$$
\mathcal{GD}^m (A) = \mathcal{GD}^{m+1}(A) = \mathcal{GD}^{m+2}(A) = \cdots 
$$
for any natural number $n \in \N$ and any subanalytic set $A \subset \R^n$ 
with $0 \in \overline{A}$ ?
\end{ques}

In \S 4 we shall give a negative example to Question \ref{SC}.

\bigskip

\section{Preliminaries}\label{prelim}

\subsection{Geometric Directional Bundle}\label{GDB}
We first recall the direction set.

\begin{defn}\label{direction}
Let $A$ be a subset of $\R^n$, 
and let $p \in \R^n$ such that $p \in \overline{A}$, 
where $\overline{A}$ denotes the closure of $A$ in $\R^n$.
We define the {\em direction set} $D_p(A)$ of $A$ at $p$ by
$$
D_p(A) := \{a \in S^{n-1} \ | \
\exists  \{ x_i \} \subset A \setminus \{ p \} ,
\ x_i \to p \in \R^n  \ \text{s.t.} \
{x_i - p \over \| x_i - p \| } \to a, \ i \to \infty \}.
$$
Here $S^{n-1}$ denotes the unit sphere centred at $0 \in \R^n$.

For a subset $D_p(A) \subset S^{n-1}$, we denote by $LD_p(A)$
the half-cone of $D_p(A)$ with $0 \in \R^n$ as the vertex, and call it the {\em real tangent cone} of $A$ at $p$:
$$
LD_p(A) := \{ t a \in \R^n\ | \ a \in D_p(A), \ t \ge 0 \}.
$$
For $p \in \overline{A}$, for simplicity, we put $L_pD(A) :=p+ LD_p(A)$, and call it the {\em geometric tangent cone} of $A$ 
at $p$.

In the case where $p = 0 \in \R^n$, we write
$D(A) := D_0(A)$ and $LD(A) := L_0D(A)$ for short.
\end{defn}

We next introduce the notion of geometric directional bundle.

\begin{defn}\label{gdirectionset} 
Let $A \subset \R^n$, and let $W \subset \R^n$ such that 
$\emptyset \ne W \subset \overline{A}$. 
We define the {\em direction set} $D_W(A)$ of $A$ over $W$ by
$$
D_W (A) := \bigcup_{p \in W} (p, D_p(A)) \subseteq W \times S^{n-1}\subseteq \R^n\times S^{n-1}.
$$

We call 
$$\mathcal{GD}_W(A) := W\times S^{n-1}\cap \overline{D_A(A)},$$
the {\em geometric directional bundle} of $A$ over $W,$
where $\overline{D_A(A)}$ denotes the closure  of $D_A(A)$
in $\R^n \times S^{n-1}$.

Let $\Pi : \R^n \times S^{n-1} \to S^{n-1}$ be the canonical
projections defined by $\Pi(x,a)=a$.
In the case where $W = \{ p \}$ for $p \in \overline{A}$, we write
$$\mathcal{GD}_p(A) := \Pi(\mathcal{GD}_{\{p\}}(A))\subseteq S^{n-1}.$$  
This is the set of all possible limits of directional sets $D_q(A), q\in A, q\to p.$


We consider the half cone of $\mathcal{GD}_p (A)$ with $0\in \R^n$ as the vertex
$$
L\mathcal{GD}_p(A) :=  \{ t v \in \R^n \ | \ v \in \mathcal{GD}_p(A), \ t \ge 0 \}.
$$ We call it the {\em real tangent bundle cone} of $A$ at $p$.

Similarly we define the {\em geometric bundle cone} of $A$ at $p,$ to be its translation by $p$
$$
L_p\mathcal{GD}(A) := p+ L\mathcal{GD}_p(A) .
$$ 
In the case where $p = 0 \in \R^n$, we simply write
$\mathcal{GD}(A) := \mathcal{GD}_0(A)$, and 
$L\mathcal{GD}(A)  := L_0\mathcal{GD}(A)$.
\end{defn}

Let $A \subset \R^n$ such that $0 \in \overline{A}$.
For $m = 2, 3, \cdots$, we define
$$
\mathcal{GD}^m(A) := \mathcal{GD}(L\mathcal{GD}^{m-1}(A)).
$$

\begin{rem} Let $A \subset \R^n$, and let $p \in \R^n$ such that $p\in \overline A$.  
If $A$ is subanalytic, then $D_p(A)\subset \mathcal{GD}_p(A).$ 
In particular, in this case we have
$$
D(A)\subseteq \mathcal{GD}^m(A)\subseteq \mathcal{GD}^{m+1}(A) 
$$
for any $m \in \N$.
\end{rem}
\subsection{Basic properties}\label{basic} 
Let $A \subset \R^n$ be a subanalytic set such that 
$0 \in \overline{A}$. 

\vspace{2mm}

{\bf (I)} $L\mathcal{GD}(A)$ is a subanalytic subest of $\R^n$.

\vspace{2mm}

{\bf (II)} The following equality holds:
$$
L\mathcal{GD}(A \setminus \{ 0\}) = L\mathcal{GD}(A).
$$

{\bf (III)} Since $\mathcal{GD}(A)$ ($= \mathcal{GD}_0(A)$) is a notion 
defined locally around $0 \in \R^n$, it suffices to consider $A \cap U$ 
for a sufficiently small neighbourhood $U$ of $0 \in \R^n$ 
when we treat Question \ref{question11}.
Throughout this paper, let us denote by $B_r(0)$ an $r$-neighbourhood of $0$ in $\R^n$ for $r > 0$, 
that is an open ball of radius $r$ with $0 \in \R^n$ as the centre.

Let us assume that $\dim A \ge 1$, and express $A \setminus \{ 0 \}$ 
in a sufficiently small $\epsilon$-neighbourhood $B_{\epsilon}(0)$ of $0 \in \R^n$ as follows: 
$$
(A \setminus \{ 0 \} ) \cap B_{\epsilon}(0) = \bigcup_{i = 1}^m C_i ,
$$
where $C_i$ is a connected component of $(A \setminus \{ 0 \} ) \cap B_{\epsilon}(0)$ 
such that $0 \in \overline{C_i}$ ($1 \le i \le m$).
Then we have 
$$
L\mathcal{GD}(A \setminus \{ 0\}) = 
\bigcup_{i = 1}^m L\mathcal{GD}(C_i ).
$$

\begin{example}\label{INR3}

(1) Let $A_+ = \{ (x,y,z) \in \R^3 \ | \ z^2 = x^2 + y^2 \ \text{and} \ 
z > 0 \}$, and 

\noindent
$A_- = \{ (x,y,z) \in \R^3 \ | \ z^2 = x^2 + y^2 \ \text{and} \ z < 0 \}$.
Set $A := A_+ \cup A_-$.

Since $L\mathcal{GD}(A_+) = L\mathcal{GD}(A_-) = 
\{ (x,y,z) \in \R^3 \ | \ z^2 \le x^2 + y^2 \}$, we have 
$L\mathcal{GD}(A) = \{ (x,y,z) \in \R^3 \ | \ z^2 \le x^2 + y^2 \}$. 
Therefore $L\mathcal{GD}^2(A) = \R^3$. 
It follows that $A \subset \R^3$ is stabilised at degree $2$.

(2) Let $A_+ = \{ (x,y,z) \in \R^3 \ | \ z^2 = 2(x^2 + y^2) \ \text{and} \ 
z > 0 \}$, and 

\noindent
$A_- = \{ (x,y,z) \in \R^3 \ | \ y^2 = 2(x^2 + z^2) \ \text{and} \ y > 0 \}$.
Set $A := A_+ \cup A_-$.

Since $L\mathcal{GD}(A) = \R^3$, 
$A \subset \R^3$ is stabilised at degree $1$.
\end{example}

{\bf (IV)}  Let $\dim A = k$, $0 \le k \le n$. 
A subanalytic set admits a locally finite Whitney stratification (H. Hironaka \cite{hironaka}). 
(See H. Whitney \cite{whitney1, whitney2} for the Whitney stratification.) 
Therefore, locally at $0 \in \R^n$ $A$ admits a finite Whitney stratification $\mathcal{S}(A)$ 
such that for any stratum $Q_i \in \mathcal{S}(A)$, $0 \in \overline{Q_i}$.
For this finite stratification $\mathcal{S}(A)$, we have the following properties:

\vspace{2mm}

(i) If $Q_i \subset \overline{Q_j}$ for $Q_i, Q_j \in \mathcal{S}(A)$, then 
$L\mathcal{GD}(Q_i) \subset L\mathcal{GD}(Q_j)$.

\vspace{2mm}

(ii) 
Let $\Omega:= \{ Q_i \in \mathcal{S}(A) | \nexists \, Q_j \ (j \ne i) 
\ \text{s.t.} \ Q_i \subset \overline Q_j\}$. 
Then we have 
$$
L\mathcal{GD}(A) = \bigcup L\mathcal{GD}(Q_i),
$$ 
where the union is taken over $\Omega$.


\subsection{A notation}

In this subsection we prepare for a notation which we will be often used in the sequel.

\begin{notation}\label{Lsubspace}
Let $m$ be a positive integer such that $0 < m < n$, and let $C \subset \R^n$ be a cone 
with $0 \in \R^n$ as a vertex.
If $C$ is contained in a finite union of $m$-dimensional linear subspaces of $\R^n$ 
(respectively If $C$ is not contained in any finite union of $m$-dimensional linear subspaces 
of $\R^n$), we write
$$
C \subseteq \bigcup_{\text{finite}} \R^m \ \ \ \ \
\text{(respectively}  \ \ C \nsubseteq \bigcup_{\text{finite}} \R^m ).
$$
\end{notation}

We give an example of an algebraic curve in $\R^n$ which is not contained in 
any finite union of hyperplanes. 

\begin{example}\label{algebraiccurve}
Let $\gamma : \R \to \R^n$, $n \ge 2$, be an algebraic curve defined 
by $\gamma (t) = (t, t^2, \cdots, t^n)$, 
and let $A \subset \R^n$ be the image of $\gamma$. 
Then even if we take a positive number $\epsilon$ arbitrarily small, 
$A \cap B_{\epsilon}(0)$ is not contained in any finite union of 
hyperplanes $\R^{n-1} \subset \R^n$.
In fact, suppose that $A \cap B_{\epsilon}(0)$ is contained 
in a finite union of hyperplanes $\bigcup_i H_i$.
Here, we may assume that each $H_i$ includes $0 \in \R^n$ 
and infinitely many points of $A \cap B_{\epsilon}(0)$.
Let one of such hyperplanes $H_i$ be given by 
$$
C_1 x_1 + C_2 x_2 + \cdots + C_n x_n = 0 \ \ \text{for} \ \
(C_1, C_2, \cdots , C_n) \ne (0, 0, \cdots, 0).
$$
By assumption, the equation $C_1 t + C_2 t^2 + \cdots + C_n t^n = 0$ 
has infinitely many roots.
It follows that $C_1 = C_2 = \cdots = C_n = 0$. 
This is a contradiction.
Thus we have
$$
A \cap B_{\epsilon}(0) \nsubseteq \bigcup_{\text{finite}} \R^{n-1}. 
$$

On other hand, we have $L\mathcal{GD}^k(A) = \R \subset \R^n$ for 
any $k \in \N$. 
Note that $\mathcal{GD}$ is stabilised at  degree $1$ for  $A$, 
and $\dim A = \dim L\mathcal{GD}(A) = 1$. 
\end{example}

\bigskip

\section{Proof of Theorem \ref{weakconj}}


\subsection{Decomposition of a subanalytic set}
In this subsection we describe  a finite decomposition of  the subanalytic sets
into subanalytic subsets equipped with some dimensional conditions.

We first show an important property on $\mathcal{GD}$.

\begin{prop} Let $A, \ B \subset \R^n$ be subanalytic sets such that 
$0 \in \overline{A} \cap \overline{B}$. 
Then the following equality holds.
$$
\mathcal{GD}(A \cup B) = \mathcal{GD}(A) \cup \mathcal{GD}(B)
$$ 
\end{prop}

\begin{proof}
One inclusion is trivial hence we concentrate in showing that $\mathcal{GD}(A \cup B) \subseteq \mathcal{GD}(A) \cup \mathcal{GD}(B)$.

In fact we show that $\overline{D_A(A)}\cup \overline{D_B(B)}=\overline{D_{A\cup B}(A\cup B)}$ so we will suffice to show 
$$
\overline{D_{A\cup B}(A\cup B)} \subseteq \overline{D_A(A)}\cup \overline{D_B(B)}=\overline{D_A(A)\cup D_B(B)}.
$$

Let  $x\in A\cup B$,  $D_x(A\cup B)=D_x(A)\cup D_x(B)$. If $x\in A\cap B$ we clearly have that 
$D_x(A\cup B)=D_x(A)\cup D_x(B)\subseteq D_{A\cup B}(A\cup B)$ and we are done.

The problem is that in general we may have $x\in A\cup B$  but $x \in \overline B \setminus B$ for instance. 
Then $D_x(A\cup B)=D_x(A)\cup D_x(B)$ but $D_x(B)$ does not count in $\mathcal{GD}(B)$. 

However in the subanalytic case $D_x(B)$ is approximated by $D_y(B)$ as $y\to x, y\in B$  i.e. 
$D_x(B) \subseteq \overline {D_B(B)}$ and thus we proved the desired inclusion.
\end{proof}

From the definition of $\mathcal{GD}^q$, we can show the following corollary 
of the above proposition by induction on $q$.

\begin{cor} 
Let $A, \ B \subset \R^n$ be subanalytic sets such that 
$0 \in \overline{A} \cap \overline{B}$, and let $q$ be a positive integer.
Then the following equality holds.
$$
\mathcal{GD}^q(A \cup B) = \mathcal{GD}^q(A) \cup \mathcal{GD}^q(B)
$$ 
\end{cor}

We can show a more general statement of the above corollary by induction 
on the number of subanalytic sets.

\begin{cor}\label{finiteunion} Let $A_1, \cdots , A_s \subset \R^n$ 
be subanalytic sets such that $0 \in \overline{A_1} \cap \cdots \cap \overline{A_s}$, 
and let $q$ be a positive integer.
Then the following equality holds.
$$
\mathcal{GD}^q(A_1 \cup \cdots \cup A_s ) 
= \mathcal{GD}^q(A_1 ) \cup \cdots \cup \mathcal{GD}^q(A_s )
$$ 
\end{cor}

Let $A \subset \R^n$ be a subanalytic set, and let $p \in A$.
Since $A$ admits a locally finite Whitney stratification with analytic strata, 
there exist a positive number $\epsilon_0 > 0$ and an integer $k$ 
with $0 \le k \le n$ such that for any $0 < \epsilon < \epsilon_0$, 
$A \cap B_{\epsilon}(p)$ is a $k$-dimensional subanalytic set, 
where $B_{\epsilon}(p)$ denotes an $\epsilon$-neighbourhood of $p$ in $\R^n$. 
We call this k {\em the local dimension of $A$ at $p$}. 
Using the local dimension, we next prepare for a dimensional condition.

\begin{defn}
Let $A \subset \R^n$ be a $k$-dimensional subanalytic set ($0 \le k \le n$). 
We say that $A$  is {\em genuinely $k$-dimensional}, 
if for any $p \in A$, the local dimension of $A$ at $p$ is $k$. 
\end{defn}

\begin{example}\label{Whitneyumbrella}
Let $A \subset \R^3$ be the Whitney umbrella. 
Then $A$ is a $2$-dimensional algebraic set.
Since the local dimension at  the points on the handle of the umbrella is $1$, 
the Whitney umbrella $A$ is not genuinely $2$-dimensional.
\end{example}

\begin{lem}\label{elementarylemma}
Let $A \subset \R^n$ be a {\em genuinely $k$-dimensional}  subanalytic set, $0 \le k \le n$,   and let $R_k$ be the set of regular points of $A$.
Then $L\mathcal{GD}(A)$ can be expressed as a union of $k$-planes, 
and $L\mathcal{GD}(A) = L\mathcal{GD}(R_k)$.
\end{lem}

\begin{proof}
Since $A$ is genuinely $k$-dimensional, the closure $\hat{R}_k$ of $R_k$ 
in $A$ coincides with $A$. 
As mentioned in Property (IV), locally at $0 \in \R^n$ $A$ admits a finite 
Whitney stratification $\mathcal{S}(A)$ such that for any stratum 
$Q_i \in \mathcal{S}(A)$, $0 \in \overline{Q_i}$, and we may assume that 
$\mathcal{S}(A)$ is compatible with $R_k$ and $A \setminus R_k$.
Since $A$ is genuinely $k$-dimensionnal, $\mathcal{S}(A)$ satisfies 
the assumption of (ii) in Property (IV).
Therefore we have 
$$
L\mathcal{GD}(A) = \bigcup L\mathcal{GD}(Q_j),
$$ 
where the union is taken over 
$\Omega_k := \{ Q_i \in \mathcal{S}(A) | \dim Q_i = k\}$. 
It follows that $L\mathcal{GD}(A)$ can be expressed as a union of 
$k$-planes. 
In addition, since 
$$
\bigcup_{Q_j \in \Omega_k} Q_j \subset R_k ,
$$ 
we have $L\mathcal{GD}(A) = L\mathcal{GD}(R_k)$. 
\end{proof}

Let $A \subset \R^n$ be a $k$-dimensional subanalytic set ($0 \le k \le n$) 
such that $0 \in \overline{A}$. 
Define the following subsets of $A$,
$$
A_i := \{ p \in A : \text{the local dimension of} \ A \ 
\text{at} \ p \ \text{is} \ i \}
$$
for $0 \le i \le k$. 
As above, locally at $0 \in \R^n$ $A$ admits a finite Whitney stratification $\mathcal{S}(A)$. 
We are assuming the frontier condition (see J. N. Mather \cite{mather}) 
for the Whitney stratification.
Therefore the local dimension is constant on each stratum of $\mathcal{S}(A)$.
It follows that if $A_i$, $0 \le i \le k$, is non-empty, $A_i$ is a genuinely  $i$-dimensional 
subanalytic set.
In particular, since $A$ is a $k$-dimensional subanalytic set, $A_k$ is non-empty and 
a genuinely $k$-dimensional subanalytic set.

Set
$$
\Lambda_0 (A) := \{ \ 0 \le i \le k \ | \ 0 \in \overline{A_i} \} .
$$
Then we can locally express $A$ as a disjoint union  around $0 \in \R^n$ as follows:
\begin{equation}\label{localunion}
A \cap B_{\epsilon}(0) 
= \bigcup_{i \in \Lambda_0 (A)} A_i \cap B_{\epsilon}(0) 
\end{equation}
for a sufficiently small $\epsilon > 0$. 
Since $0 \in \overline{A}$, $\Lambda_0 (A)$ is non-empty, thus   we may  consider  $m_0 (A) := \min \Lambda_0 (A)$. 

\begin{rem}\label{remark30}
If $0 \notin A$, $A_0$ is empty locally at $0 \in \R^n$. 
Therefore $0 \notin \Lambda_0 (A)$. 

In addition, if $0 \in A_i$ for some positive integer $i$, 
then $0 \notin \Lambda_0 (A)$. This shows that $m_0(A)=0$ only if $0$ is an isolated point of $A$.
\end{rem}

Let $q$ be a positive integer.
By (II) and Corollary \ref{finiteunion}, we have the following formula.
\begin{equation}\label{union}
L\mathcal{GD}^q (A)
= L\mathcal{GD}^q (A \setminus \{ 0\})
= \bigcup_{i \in \Lambda_0 (A)} L\mathcal{GD}^q(A_i \setminus \{ 0\})
= \bigcup_{i \in \Lambda_0 (A)} L\mathcal{GD}^q(A_i).
\end{equation}


\subsection{Decomposition of a geometric bundle cone}

In this subsection we discuss some properties concerning  a specific decomposition of  the geometric bundle cones.

\vspace{2mm}

{\bf (V)} Let $A \subset \R^n$ be a $k$-dimensional subanalytic set, $0 \le k \le n, $ 
such that $0 \in \overline{A}$.
Write $A$ as in the equality (\ref{localunion}) above and denote by $R_i$ the subset of all regular points  in $A_i, i\in \Lambda_0 (A)$. 
$$
L\mathcal{GD}(A)=\bigcup_{i \in \Lambda_0 (A)} L\mathcal{GD}(R_{i}) 
$$ 
and each $L\mathcal{GD}(R_{i}) $ is a union of $i-$planes.

\vspace{2mm}

{\bf (VI)} Let $A$ be a genuinely $k-$dimensional subanalytic set, $0 \le k \le n$, 
such that $0 \in \overline{A}$, and let $R_k$ be the set of regular points of $A$.
By Lemma \ref{elementarylemma}, $B:=L\mathcal{GD}(A) =  L\mathcal{GD}(R_{k})$ 
and $B$ is a union of $k-$planes $B=\cup_{i\in I}\pi_i$. 
In $B$ we distinguish two kind of points:
$$
B_+:=\{x\in B | \forall  V \text{open in $B$} , x\in V, V\cap \pi_i\neq\emptyset ,\,  \text{for infinitely many indices}\, i\in I\}, 
$$ 
$$
B_0:=\{x\in B | \exists V \text{open in $B$} , x\in V, V\cap \pi_i\neq\emptyset ,\,  \text{only for finitely many indices}\, i\in I\},
$$ 
thus 
$B= B_+\cup B_0$. By definition $B_0$ is open in $B$ and thus its complement $B_+$ is closed. Note that, as  $B=\cup_{i\in I}\pi_i$,  $B_0$ has to be contained in a finite union of  $k-$planes $\pi_i$ from $B$ (it is a subset of a union of $k-$planes and the local dimension is $k$). As such, $L\mathcal{GD}(B)=L\mathcal{GD}(B_+)\cup L\mathcal{GD}(B_0)$ where $L\mathcal{GD}(B_0)$, if not empty, is a finite union of $k-$planes, thus of dimension $k$. In fact we know that $B_0$ is genuinely $k-$dimensional by construction so $L\mathcal{GD}(B_0)=L\mathcal{GD}(R), \, R$  its regular points, which has to be a finite union of connected components and each component has to be a subset of a participating plane in $B_0$, so $L\mathcal{GD}(R)$ must be the finite union of those planes. Now if $B_+$ is not empty it has to have dimension bigger than $k$ and thus we can apply the decomposition above to conclude that 
$L\mathcal{GD}(B_+)=\cup L\mathcal{GD}(R_j)$ where $R_j$ are $j-$dimensional analytic manifolds with $j>k$.

\begin{example}\label{GDBpartition}
Let $A \subset \R^3$ be a 2-dimensional algebraic set defined by
$$
A := \{ (x,y,z) \in \R^3 \ | \  xy(x^2 + y^2 -z^2) = 0 \} .
$$
Then $A$ is a genuinely $2$-dimensional subanalytic set, and we have 

\vspace{3mm}

\qquad $B_0 = \{ (x,y,x) \in \R^3 \  | \ xy = 0 \ \& \ z^2 > x^2 + y^2 \}$,

\vspace{2mm}

\qquad $B_{+}  = \{ (x,y,z) \in \R^3 \ | \ z^2 \le x^2 + y^2 \} .$

\vspace{2mm}

Let $R_2$ be the set of regular points of $A$, and let $B := L\mathcal{GD}(A) = L\mathcal{GD}(R_2)$. 
$B = B_0 \cup B_+$ and $B$ is a union of 2-planes.
Note that $B_0$ is contained in a union of two $2$-planes.

$B_0$ is $2$-dimensional  and open in $B$. 
On the other hand, $B_{+}$ is $3$-dimensional and closed in $B$. 
In addition, $L\mathcal{GD}(B_0) = \{ (x,y,z) \in \R^3 \ | \ xy = 0\}$  
is a union of two $2$-planes.
\end{example}


\subsection{Observations}\label{observation}
In this subsection we observe that Theorem \ref{weakconj} holds in the case where $n \le 3$.
We first show the case $n = 1$. 

\begin{ass}\label{n=1}
If $n = 1$, $\mathcal{GD}$ is stabilised at degree $1$.
\end{ass}

\begin{proof}
Let $A \subset \R$ be a subanalytic set such that $0 \in \overline{A}$. 
Then $A$ is $\{ 0\}$, $[0, \epsilon )$, $(0, \epsilon )$, $(-\epsilon ,0]$, 
$(-\epsilon ,0)$ ($\epsilon$ is a sufficiently small positive number), 
$(-\epsilon_1 ,\epsilon_2)$ or $(-\epsilon_1 ,\epsilon_2) \setminus \{ 0 \}$ 
($\epsilon_1$ and $\epsilon_2$ are sufficiently small positive numbers) 
in a small neighbourhood of $0 \in \R$.
In case $A = \{ 0\}$, $L\mathcal{GD}(A) = \emptyset$.
Otherwise, we have $L\mathcal{GD}(A) = \R$.
It follows that $\mathcal{GD}$ is stabilised at degree $1$.
\end{proof}

We next show the case $n = 2$. 

\begin{lem}\label{onedimensional}
Let $n \ge 2$.
Suppose that $\dim A = 1$.
Then $\mathcal{GD}$ is stabilised at degree $1$.
\end{lem}

\begin{proof}
Let us express $A \setminus \{ 0 \}$ in a sufficiently small 
$\epsilon$-neighbourhood $B_{\epsilon}(0)$ of $0 \in \R^n$ like 
in (2.2)(III):
$$
(A \setminus \{ 0 \} ) \cap B_{\epsilon}(0) = \bigcup_{i = 1}^m C_i ,
$$
where $C_i$ is a connected component of $(A \setminus \{ 0 \} ) \cap B_{\epsilon}(0)$ 
such that $0 \in \overline{C_i}$ ($1 \le i \le m$).
Here $C_i$ is a subanalytic curve.
Therefore each $L\mathcal{GD}(C_i) = \R \subset \R^n$ ($1 \le i \le m$) 
in our notation.
(Some $\R$'s may coincide.)
It follows that $\mathcal{GD}$ is stabilised at degree $1$.
\end{proof}

\begin{ass}\label{n=2}
If $n = 2$, $\mathcal{GD}$ is stabilised at degree $1$.
\end{ass}

\begin{proof}
In the case where $\dim A = 0$, $L\mathcal{GD}(A) = \emptyset$.
In the case where $\dim A = 2$, $L\mathcal{GD}(A) = \R^2$.
In these cases it is obvious that $\mathcal{GD}$ is stabilised at degree $1$.

By Lemma \ref{onedimensional}, 
$\mathcal{GD}$ is stabilised at degree $1$ in the case where $\dim A = 1$.
\end{proof}

\begin{ass}\label{n=3}
If $n = 3$, $\mathcal{GD}$ is stabilised at degree $2$. 
\end{ass}

\begin{proof}
In the cases where $\dim A = 0$ and $\dim A = 1$,
$\mathcal{GD}$ is stabilised at degree $1$ as seen in the proof of Assertion \ref{n=2}. 
In the case where $\dim A = 3$, $L\mathcal{GD}(A) = \R^3$. 
Therefore $\mathcal{GD}$ is stabilised at degree $1$.

Let us consider the case $\dim A = 2$.
As in Lemma \ref{onedimensional}, let us express $A \setminus \{ 0 \}$ 
in a sufficiently small $\epsilon$-neighbourhood $B_{\epsilon}(0)$ of $0 \in \R^3$:
$$
(A \setminus \{ 0 \} ) \cap B_{\epsilon}(0) = \bigcup_{i = 1}^m C_i .
$$
We may assume that $C_1, \cdots , C_k$, $1 \le k \le m$, are genuinely $2$-dimensional,
and $C_{k+1}, \cdots , C_m$ are $1$-dimensional. 
In the case where $k+1 \le i \le m$, namely in the $1$-dimensional case,  
we have $L\mathcal{GD}(C_i) = \R \subset \R^3$ in our notation.
Therefore $\mathcal{GD}$ is stabilised at degree $1$.

We next consider the $2$-dimensional case. 
By Lemma \ref{elementarylemma}, for each $1 \le i \le k$ $L\mathcal{GD}(C_i)$ is 
a union of $2$-planes.
Let 
$$
L\mathcal{GD}(C_i) \subseteq \bigcup_{\text{finite}} \R^2, \ 1 \le i \le s, \ \text{and} \ 
L\mathcal{GD}(C_i) \nsubseteq \bigcup_{\text{finite}} \R^2, \ s + 1 \le i \le k.
$$
In the case where $1 \le i \le s$, we have $L\mathcal{GD}(C_i) = \bigcup_{\text{finite}} \R^2$. 
Therefore $\mathcal{GD}$ is stabilised at degree $1$.
On the other hand, in the case where $s + 1 \le i \le k$, the subanalytic set 
$L\mathcal{GD}(C_i)$ contains infinitely many $2$-planes.
Therefore we have $\dim L\mathcal{GD}(C_i) = 3$, and 
$L\mathcal{GD}^2(C_i) = \R^3$.
Thus $\mathcal{GD}$ is stabilised at degree $2$.

It follows that $\mathcal{GD}$ is stabilised at degree $2$ in the case of $n = 3$.

\end{proof}


\subsection{Proof in the general case}\label{Wconjecture}
We first recall some fundamental properties, seen in \ref{observation}.

\begin{pro}\label{property1} 
Let $A \subset \R^n$ be a subanalytic set such that $0 \in \overline{A}$.  
Then we have the following.

(1) If $\dim A = 0$, $\mathcal{GD}$ is stabilised at degree $1$. 

(2) If $\dim A = 1$, $\mathcal{GD}$ is stabilised at degree $1$. 

(3) If $\dim A = n$, $\mathcal{GD}$ is stabilised at degree $1$. 
\end{pro}

Before we start the proof of Theorem \ref{weakconj}, we prepare for 
some lemmas. 
The following lemma is an immediate consequence of Lemma \ref{elementarylemma}.

\begin{lem}\label{lemmap1}
Let $A \subset \R^n$ be a genuinely $k$-dimensional subanalytic set, 
$1 \le k \le n$, such that $0 \in \overline{A}$. 
Then we have $m_0 (L\mathcal{GD}(A)) \ge k$. 
\end{lem}

\begin{lem}\label{lemmap2}
Let $A \subset \R^n$ be a genuinely $k$-dimensional subanalytic set, 
$1 \le k \le n$, such that $0 \in \overline{A}$. 
Suppose that $\dim L\mathcal{GD}(A) = k$. 
Then $L\mathcal{GD}(A)$ can be expressed as a union of some finite $k$-planes.
\end{lem}

\begin{proof}
By Lemma \ref{elementarylemma}, $L\mathcal{GD}(A)$ can be expressed 
as a union of $k$-planes.
Let us assume that the subanalytic set $L\mathcal{GD}(A)$ cannot be expressed 
as a union of finite $k$-planes. 
Then it is a union of infinitely many different $k$-planes.
It implies that $\dim L\mathcal{GD}(A) > k$. 
This contradicts the assumption $\dim L\mathcal{GD}(A) = k$. 
Thus $L\mathcal{GD}(A)$ can be expressed as a union of some finite $k$-planes.
\end{proof}


Related to the above lemma, we can present a proposition. 
Let $A \subset \R^n$ be a subanalytic set such that $0 \in \overline{A}$. 
Like (\ref{localunion}), let us locally express $A$ around $0 \in \R^n$ as follows:
$$
A \cap B_{\epsilon}(0) 
= \bigcup_{i \in \Lambda_0 (A)} A_i \cap B_{\epsilon}(0) 
$$
for a sufficiently small $\epsilon > 0$.  
Then we have the following proposition.

\begin{prop}\label{dimLGD}
Let $i_0$ be the smallest integer in $\Lambda_0 (A)$ such that 
$L\mathcal{GD}(A_{i_0})$ is not contained in any finite union of hyperplanes.
Then we have $\dim L\mathcal{GD}^{n-i_0} (A) = n$. 
\end {prop}

We next introduce some notations  concerning $L\mathcal{GD}(A)$.

\begin{notation}  
Let $A \subset \R^n$ be a subanalytic set such that $0 \in \overline{A}$, 
and let $m_0(A) = k$, $1 \le k < n$. 
Then let us denote by $L\mathcal{GD}(A)_{= k}$ (respectively 
$L\mathcal{GD}(A)_{\ge k+1}$) the set of points of $L\mathcal{GD}(A)$ 
at which the local dimension of it is $k$ (respectively bigger than or equal 
to $k + 1$).
\end{notation}

\begin{rem}\label{remark31}
Let $A \subset \R^n$ be a subanalytic set such that $0 \in \overline{A}$, 
and let $m_0(A) = k$, $1 \le k < n$. 
By Lemma \ref{lemmap1}, we have
$$
L\mathcal{GD}(A) = L\mathcal{GD}(A)_{= k} \cup L\mathcal{GD}(A)_{\ge k + 1}.
$$
\end{rem}

\begin{example}\label{reversal}
Let $A \subset \R^6$ be a 5-dimensional algebraic set defined by
$$
A := \{ (x,y,z,u,v,w) \in \R^6 \ | \ 
x \{ (x^2 + y^2 -z^4)^2 + (u^2 + v^2 -w^4 )^2 \} = 0 \} .
$$
Let us denote by $A_{=5}$ (respectively $A_{=4}$) the set of points of $A$ 
at which the local dimension of it is $5$ (respectively $4$) as above.
Then we can express $A$ as the union of $A_{=5}$ and $A_{=4}$, where   

\vspace{3mm}

$A_{=5} = \{ (x,y,z,u,v,w) \in \R^6 \ | \ x = 0 \}$, 

\vspace{3mm}

$A_{=4} = \{ (x,y,z,u,v,w) \in \R^6 \ | \ 
(x^2 + y^2 -z^4)^2 + (u^2 + v^2 -w^4 )^2 = 0 \}$

\vspace{3mm}

\qquad $\setminus \{ (x,y,z,u,v,w) \in \R^6 \ | \ x = 0 \ \& \  
(x^2 + y^2 -z^4)^2 + (u^2 + v^2 -w^4 )^2 = 0 \} $.

\vspace{3mm}

We can see that $A_{=5} = \R^5 \subset \R^6$ and 
$L\mathcal{GD}(A_{=4}) = \R^6$. 
Therefore $\mathcal{GD}$ is stabilised at degree $0$ for $A_{=5}$ in some sense, 
and stabilised at degree $1$ for $A_{=4}$. 
It follows that $\mathcal{GD}$ is stabilised at degree $1$ for $A$. 
We note that the lower dimensional subset $A_{=4}$ of $A$ takes 
an essential role in this stabilisation for $A$. 
\end{example}

We prepare for one more lemma.

\begin{lem}\label{lemmap3}
Let $A \subset \R^n$ be a genuinely $k$-dimensional subanalytic set, 
$1 \le k \le n$, such that $0 \in \overline{A}$. 
Then $L\mathcal{GD}(A)_{=k}$ is contained in a union of some finite $k$-planes.
\end{lem}

\begin{proof}
By Lemma \ref{elementarylemma}, $B := L\mathcal{GD}(A)$ is a union of 
$k$-planes $B = \cup_{i\in I} \pi_i$. 
We consider the decomposition of $B$ in (VI), namely $B = B_0 \cup B_+$. 
On the other hand, we have another decomposition 
$B = L\mathcal{GD}(A)_{=k} \cup L\mathcal{GD}(A)_{\ge k+1}$. 
By definition, we can easily see that $B_0 \subset L\mathcal{GD}(A)_{=k}$. 
Pick a point $P \in L\mathcal{GD}(A)_{=k}$. 
Then there exists a subanalytic neighbourhood $W \subset \R^n$ of $P$ such that 
the $k$-dimensional subanalytic set $W \cap B$ admits a finite Whitney stratification 
$\mathcal{S}_P(B)$ satisfying the following conditions:

(1)  Any $k$-dimensional stratum of $\mathcal{S}_P (B)$ is analytically diffeomorphic 
to a $k$-plane.

(2) For any $k$-dimensional stratum $\sigma \in \mathcal{S}_P (B)$, $P$ is 
 in the closure of $\sigma$.
 
Suppose that $P \in B_+$. 
Then there exists an infinite subset $I_1$ of $I$ such that
$W \cap \pi_i \ne \emptyset$ for $i \in I_1$ and  
$W \cap B = W \cap \cup_{i \in I_1} \pi_i$ ($= \cup_{i \in I_1} (W \cap \pi_i )$). 
Note that $\cup_{i \in I_1}\pi_i$ cannot yield a point in $W$ at which the local dimension 
of $B$ is bigger than $k$.
Then the union of strata of $\mathcal{S}_P (B)$ does not coincide with 
$\cup_{i \in I_1} (W \cap \pi_i )$.
Therefore $P \in B_0$. 
It follows that $L\mathcal {GD}(A)_{=k} = B_0$. 
Thus, by (VI), $L\mathcal{GD}(A)_{=k}$ is contained 
in a union of some finite $k$-planes.
\end{proof}

Now let us start the proof of Theorem \ref{weakconj}.  
By Assertions \ref{n=1} and \ref{n=2}, it suffices to show the theorem
in the case where $n \ge 3$. 
After this, we assume $n \ge 3$. 
Under this assumption, we have the following result.

\begin{prop}\label{main}
Let $A \subset \R^n$ be a subanalytic set such that $0 \in \overline{A}$. 
Suppose that $2 \le m_0 (A) \le n$. 
Then $\mathcal{GD}$ is stabilised at degree $n - (m_0 (A) - 1)$.
\end{prop}
 
\begin{proof}
Let us show this proposition by downward induction on $m_0 (A)$. 
By Property \ref{property1} (3), $\mathcal{GD}$ is stabilised at degree $1$ 
for any genuinely $n$-dimensional subanalytic set $A$ with $0 \in \overline{A}$. 
Therefore the proposition holds in the case  where $m_0 (A) = n$.
 
Let $i$ be a positive integer with $2 \le i \le n - 1$.
Let us assume that this proposition holds for any subanalytic set $B$ 
with $0 \in \overline{B}$ such that $m_0 (B) \ge i + 1$.
Let $A \subset \R^n$ be an arbitrary subanalytic set with $m_0 (A) = i$
such that $0 \in \overline{A}$.  

We first consider the case where $L\mathcal{GD}(A)_{=i} = \emptyset$.
Then $m_0(L\mathcal{GD}(A)) \ge i + 1$. 
By assumption on induction, for the subanalytic set $L\mathcal{GD}(A)$, 
$\mathcal{GD}$ is stabilised at degree $n - (m_0(L\mathcal{GD}(A)) - 1)$.    
Therefore, for $L\mathcal{GD}(A)$, $\mathcal{GD}$ is stabilised 
at degree $n - i$. 
It follows that for $A$, $\mathcal{GD}$ is stabilised at degree 
$n - i + 1 = n - (i - 1) = n - (m_0 (A) - 1)$. 

We next consider the case where $L\mathcal{GD}(A)_{=i} \ne \emptyset$. 
If $L\mathcal{GD}(A)_{\ge i+1} = \emptyset$, then 
$\cup_{j \in \Lambda_0(A) \setminus \{ i\}} A_j = \emptyset$. 
Therefore $A = A_i$ is a genuinely $i$-dimensional subanalytic set.
Since $dim L\mathcal{GD}(A) = i$, by Lemma \ref{lemmap2}, 
$L\mathcal{GD}(A)$ can be expressed as a union of some finite $\R^{i}$. 
Therefore, for $A$, $\mathcal{GD}$ is stabilised at degree  $1$. 
Since 
$$
n - (m_0 (A) - 1) = n - i + 1 \ge n - (n - 1) + 1 = 2 > 1, 
$$
$\mathcal{GD}$ is stabilised at degree $n - (m_0 (A) - 1)$.

If $L\mathcal{GD}(A)_{\ge i+1} \ne \emptyset$, then by Remark \ref{remark31}, 
we have
$$
L\mathcal{GD}(A) = L\mathcal{GD}(A)_{= i} \cup L\mathcal{GD}(A)_{\ge i + 1}. 
$$
Note that this union is a disjoint one.
By Lemma \ref{lemmap1}, $L\mathcal{GD}(A)_{=i} \subset L\mathcal{GD}(A_i)$. 
It follows that $L\mathcal{GD}(A)_{=i} = L\mathcal{GD}(A_i)_{=i}$. 
Since $A_i$ is a genuinely $i$-dimensional subanalytic set,
by Lemma \ref{lemmap3}, $L\mathcal{GD}(A)_{= i}$ is contained in a union of 
some finite $\R^i$.
Therefore for $L\mathcal{GD}(A)_{= i}$, $\mathcal{GD}$ is stabilised 
at degree 1. 
On the other hand, $m_0 (L\mathcal{GD}(A)_{\ge i+1}) \ge i + 1$. 
Using the same argument as in the case of $L\mathcal{GD}(A)_{= i} = \emptyset$, 
we can see that for $L\mathcal{GD}(A)_{\ge i+1}$, $\mathcal{GD}$ is 
stabilised at degree $n - i$. 
Therefore for $L\mathcal{GD}(A)$, $\mathcal{GD}$ is stabilised at degree 
$\max \{ 1, n - i\} = n - i$. 
It follows that for $A$, $\mathcal{GD}$ is stabilised at degree 
$n - i + 1 = n - (m_0 (A) - 1)$. 
\end{proof}

\begin{proof}[Proof of Theorem \ref{weakconj}]
By Proposition \ref{main}, if $2 \le m_0 (A) \le n$, then  $\mathcal{GD}$ is stabilised 
at degree $n - (m_0 (A) - 1)$.
Since
$$
n - (m_0 (A) - 1) \le n -(2 - 1) = n - 1,
$$
$\mathcal{GD}$ is stabilised at degree $n - 1$ under the assumption that
$2 \le m_0 (A) \le n$.

We next consider tha case where $m_0 (A) = 1$.
If $A$ is a $1$-dimensional subanalytic set, then by Property \ref{property1} (2), 
$\mathcal{GD}$ is stabilised at degree $1$. 
Since $n \ge 2$, $\mathcal{GD}$ is stabilised at degree $n - 1$.
If $A$ is a subanalytic set of dimension bigger than or equal to $2$, 
we have a partition of $A$ into two subanalytic subsets $A_1$ and $A_2$, 
where $A_1$ consists of finite subanalytic curves and $A_2$ consists 
of points of $A$ at which the local dimension of it is bigger than or equal to $2$. 
Therefore for $A_1$, $\mathcal{GD}$ is stabilised at degree $1$, and 
for $A_2$, $\mathcal{GD}$ is stabilised at degree $n - 1$. 
It follows that for $A$, $\mathcal{GD}$ is stabilised at degree 
$\max \{ 1, n - 1\} = n -1$.

We lastly consider the case where $m_0 (A) = 0$. 
In this case, the dimension of $A$ locally at $0 \in \R^n$ is $0$. 
Because if the dimension of $A$ locally at $0 \in \R^n$ is bigger than $0$, 
then $0 \in \R^n$ is a point of some genuinely positive dimensional 
subanalytic subset of $A$ or $0 \notin A$.
Therefore the condition $m_0 (A) = 0$ is equivalent to the condition 
$A = \{ 0 \}$ locally at $0 \in \R^n$. 
Then, by Property \ref{property1} (1), $\mathcal{GD}$ is stabilised at degree $1$. 
It follows that $\mathcal{GD}$ is stabilised at degree $n - 1$.
\end{proof}

\bigskip

\section{Negative example to Question 2}

Let us recall Question 2.

\vspace{3mm}

\noindent
{\bf Question 2,} 
Does there exist a natural number $m \in \N$ such that 
$$
\mathcal{GD}^m (A) = \mathcal{GD}^{m+1}(A) = \mathcal{GD}^{m+2}(A) = \cdots 
$$
for any natural number $n \in \N$ and any subanalytic set $A \subset \R^n$ 
with $0 \in \overline{A}$ \ ?

\vspace{3mm}

We first prepare for some lemmas.

\begin{lem}\label{dimlemma}
Let $A \subset \R^n$ be a $k$-dimensional subanalytic set, $0 \le k \le n$,
such that $0 \in \overline{A}$. 
Then we have 
$$
k \le \dim L\mathcal{GD}(A) \le \min (2k, n) .
$$
\end{lem}

\begin{proof}
By definition, it is obvious that $k \le \dim L\mathcal{GD}(A) \le n$.
Therefore let us show the statement that $\dim L\mathcal{GD}(A) \le 2k$ if $2k < n$.

We locally express $A$ around $0 \in \R^n$ as (3.1):
$$
A \cap B_{\epsilon}(0) 
= \bigcup_{i \in \Lambda_0 (A)} A_i \cap B_{\epsilon}(0) 
$$
for a sufficiently small $\epsilon > 0$, where each $A_i$, $i \in \Lambda_0 (A)$, 
is a genuinely $i$-dimensional subanalytic set.
If we can show that $\dim L\mathcal{GD}(A_i) \le 2i$, $i \in \Lambda_0 (A)$, 
then the statement follows. 
Therefore, from the beginning, we may assume that 
$A$ is a genuinely $k$-dimensional subanalytic set.
Let $R_k$ be the set of regular points of $A$.
By Lemma \ref{elementarylemma},
$L\mathcal{GD}(A) = L\mathcal{GD}(R_k)$.
Since $\mathcal{GD}_{\{ 0\} }(R_k) = \{ 0 \} \times S^{n-1} \cap \overline{D_{R_k}(R_k)}$, 
$\dim \mathcal{GD}_{\{ 0\} }(R_k) \le 2k - 1$.
Therefore we have 
$$
\dim \mathcal{GD}(R_k) = \dim \Pi (\mathcal{GD}_{\{ 0\} }(R_k)) \le 2k - 1,
$$
where $\Pi : \R^n \times S^{n-1} \to S^{n-1}$ is the canonical projection by definition.
But in this case, since $\mathcal{GD}_{\{ 0\} }(R_k) \subset \{ 0 \} \times S^{n-1}$,
we can regard $\Pi |_{\{ 0\} \times S^{n-1}}$ as the identification map
between $\{ 0 \} \times S^{n-1}$ and $S^{n-1}$.
It follows that $\dim L\mathcal{GD}(A) = \dim L\mathcal{GD}(R_k) \le (2k - 1) + 1 = 2k$.
\end{proof}

\begin{rem}
Let $A \subset \R^n$ be a  $k$-dimensional subanalytic set. 
Then $\dim L\mathcal{GD}^2(A)\leq 4k-1$ for instance.
Indeed 
$$
L\mathcal{GD}^2(A)=L\mathcal{GD}(L\mathcal{GD}(A))=L\mathcal{GD}(Cone(B))
=Cone(L\mathcal{GD}(B)),
$$ 
where $L\mathcal{GD}(A)=Cone(B)$ and $L\mathcal{GD}(B) := \cup_{b\in B}L_b\mathcal{GD}(B)$, the union of the corresponding geometric bundles..
Now $\dim L\mathcal{GD}(A)\leq 2k$ so $\dim B \leq 2k-1$, thus by induction dim 
$L\mathcal{GD}(B)\leq 2(2k-1)$ thus finally $\dim Cone(L\mathcal{GD}(B))\leq 2(2k-1)+1=4k-1$, a bit of improvement.
Using the same observation we may get an estimation for 
$$
\dim L\mathcal{GD}^r(A)\leq 2^rk-2^{r-1}+1, \,\, r\geq 1.
$$ 
Here we used the formula $L\mathcal{GD}(Cone(B))=Cone(L\mathcal{GD}(B))$. 
Thus starting with $A,$ two dimensional in $\R^8$, such that 
$ L\mathcal{GD}(A) $ is not contained in any finite union of hyperplanes, then $L\mathcal{GD}^3(A)\neq L\mathcal{GD}^2(A)$. If we start with $A, 2$-dimensional such that $ L\mathcal{GD}(A) $ is $3$-dimensional, then
$ L\mathcal{GD}^2(A) $  has dimension $\leq 1+4=5$ thus it cannot stabilise  at $2$ in $\R^6$.
Can we find such examples in $\R^4, \R^5$?
\end{rem}

\begin{example}\label{product}
Let $A_0$ be a $2$-dimensional algebraic set in $\R^3$ defined by 
$$
A_0 := \{ (x,y,z) \in \R^3 \ | \ x^2 + y^2 = z^4 \} .
$$
Then it is easy to see that $L\mathcal{GD}(A_0 ) = \R^3$ thus it stabilises at degree $1$.

Let $A$ be a $4$-dimensional algebraic set in $\R^6$ defined by
\begin{eqnarray*}
A &=&  \{ (x,y,z,u,v,w) \in \R^6 \ | \ x^2 + y^2 = z^4 \ \& \ u^2 + v^2 = w^4 \} \\
&=& \{ (x,y,z) \in \R^3 \ | \ x^2 + y^2 = z^4 \} \times 
\{ (u,v,w) \in \R^3 \ | \ u^2 + v^2 = w^4 \} .
\end{eqnarray*}
Then we have $L\mathcal{GD}(A) = \R^3 \times \R^3 = \R^6$, thus it stabilises at degree $1$. In this way we can create examples of subspaces $A$ in any codimension, such that $ L\mathcal{GD}(A)$ is the entire space.

\end{example}

\begin{lem}\label{speciallem}
Let $A \subset \R^n$ be a $k$-dimensional subanalytic set, $0 < k < n$,
such that $0 \in \overline{A}$. 
Suppose that $L\mathcal{GD}(A)$ is not contained in any finite union of hyperplanes.
Then $L\mathcal{GD}(A) \nsubseteq A \cup F$,  where $L\mathcal{GD}(F)$ 
is contained in a finite union of hyperplanes.
\end{lem}

\begin{proof} Consider a decomposition of $A$ into a $k$-dimensional subanalytic set $A_k$ 
and a smaller dimensional subanalytic set $B$. 
Here we take a decomposition in (3,1) and $B = \cup_{i<k} A_i$. 
In this case $A_k$ is a genuinely $k$-dimensional.
If $L\mathcal{GD}(A) = L\mathcal{GD}(B) \cup L\mathcal{GD}(A_k) \subseteq A \cup F$ such that 
$L\mathcal{GD}(F)$ is contained in a finite union of hyperplanes, then there are two cases to consider.

\noindent
Case I: Let $L\mathcal{GD}(A_k)$ be a finite union of $k$-planes. 
Then $L\mathcal{GD}(B) \subset B \cup A_k \cup F$ where $L\mathcal{GD}(A_k \cup F)$ is contained
in a finite union of hyperplanes and $L\mathcal{GD}(B)$ is not contained in any finite union 
of hyperplanes, which is a contradiction by induction on dimension.

\noindent
Case II: Let $L\mathcal{GD}(A_k)$ be not contained in any finite union of hyperplanes. 
Then $L\mathcal{GD}(A_k) =A^+ \cup A^0$ such that $A^0$ is a finite union of hyperplanes 
and $A^+$ of dimension bigger than $k$ at any point. 
Accordingly $L\mathcal{GD}(A^+)$ is not contained in any finite union of hyperplanes (same as $L\mathcal{GD}(A_k)$).
We have $A^+ \subseteq A \cup F$, so $A^+ \setminus A \subseteq F$. 
Note that $A^+ \setminus A$ is dense in $A^+$ and thus
$L\mathcal{GD}(A^+) = L\mathcal{GD}(A^+ \setminus A) \subseteq L\mathcal{GD}(F)$ is contained 
in a finite union of hyperplanes by assumption, which is a contradiction again.

We conclude that $L\mathcal{GD}(A) \nsubseteq A \cup F$ where $L\mathcal{GD}(F)$ is contained in a finite union 
of hyperplanes.
\end{proof}
 
As corollaries of Lemma \ref{speciallem}, we have the following.
\begin{cor}
Let $A \subset \R^n$ be a subanalytic set
such that $0 \in \overline{A}$. Suppose that $L\mathcal{GD}(A)$ is not contained in any finite union of hyperplanes.
If $A=L\mathcal{GD}(A)$, then $A=\R^n$.
\end{cor}
\begin{cor}\label{specialcor}
Let $A \subset \R^n$ be a  subanalytic set
such that $0 \in \overline{A}$. 
Let $m$ be a positive integer such that $1 < m < n-1$.
Suppose that $L\mathcal{GD}(A)$ is not contained in any finite union of hyperplanes, and that 
$$
 L\mathcal{GD}^{m-1}(A) \neq L\mathcal{GD}^m(A) 
{\text and} \ L\mathcal{GD}^m(A) = L\mathcal{GD}^{m+1}(A).
$$
Then $L\mathcal{GD}^m(A) = \R^n.$ 
\end{cor}

\begin{example}\label{regularcurve}
Let $n$ be an a natural number with $n \ge 3$, and let $J := (-\frac{1}{2}, \frac{1}{2})$. 
We define a regular curve $\gamma : J \to S^{n-1}$ by 
$$
\gamma (t) := (t, t^2, \cdots, t^{n-1}, \sqrt{1 - (t^2 + \cdots t^{2(n-1)}}).
$$ 
Let us denote by $C$ the image of $\gamma$. 

We call the intersection of $S^{n-1}$ and a hyperplane in $\R^n$ passing through $0 \in \R^n$ 
a {\em hyperplane of $S^{n-1}$}. 
Let us show that $C \subset S^{n-1}$ is not contained in any finite union 
of hyperplanes of $S^{n-1}$.
We assume that $C$ is contained in a finite union of hyperplanes $\bigcup_i (H_i \cap S^{n-1})$ of $S^{n-1}$.
Here, each $H_i$ is a hyperplane in $\R^n$ passing through $0 \in \R^n$, and some of $H_i$'s 
include infinitely many points of $C$. 
Let one of such hyperplanes $H_i$ be given by 
$$
a_1 x_1 + a_2 x_2 + \cdots + a_n x_n = 0 \ \ \text{for} \ \
(a_1, a_2, \cdots , a_n) \ne (0, 0, \cdots, 0).
$$
By assumption, the equation 
$$
a_1 t + a_2 t^2 \cdots + a_{n-1} t^{n-1} +a_n  \sqrt{1 - (t^2 + \cdots t^{2(n-1)})} = 0
$$  
has infinitely many roots.
Therefore the polynomial equation
$$
(a_1 t + a_2 t^2 \cdots + a_{n-1} t^{n-1})^2 = (- a_n  \sqrt{1 - (t^2 + \cdots t^{2(n-1)})})^2
$$
also has infinitely many roots.
It follows that $a_n = 0$ and $a_1 = a_2 = \cdots = a_{n-1} = 0$. 
This is a contradiction.
Thus  $C \subset S^{n-1}$ is not contained in any finite union 
of hyperplanes of $S^{n-1}$.

Let $A$ be a cone of C with $0 \in \R^n$ as the vertex.
Then $\dim A = 2$ and  $\dim L\mathcal{GD}(A) = 3$. 
By construction, we can see that $L \mathcal{GD}(A) \cap S^{n-1} \supset C$. 
Therefore $L\mathcal{GD}(A)$ is not contained in any finite union of hyperplanes in $\R^n$. 

Suppose that there exists a natural number $m \in \N$ such that 
$\mathcal{GD}$ is stabilised at degree $m$ for any natural number $n \in \N$ 
and any subanalytic set $A \subset \R^n$ with $0 \in \overline{A}$. 
Then we choose a natural number $n \in \N$ so that $n > 3 \cdot 2^{m-1}$. 
We construct $A \subset \R^n$ in the above way.
By Lemma \ref{dimlemma}, we have $\dim L\mathcal{GD}^m(A) \le 3 \cdot 2^{m-1} < n$, 
which contradicts Corollary \ref{specialcor}. 
Thus there does not exist a natural number $m \in \N$ such that 
$\mathcal{GD}$ is stabilised at degree $m$ for any natural number $n \in \N$ 
and any subanalytic set $A \subset \R^n$ with $0 \in \overline{A}$. 

\end{example}

\bigskip


\end{document}